\documentclass[fleqn]{mat01}
\usepackage{times,mathtimy,amssymb,latexsym}%,krepsf,rangecite}
\begin{document}

\setcounter{page}{437}
\firstpage{437}

\def\d{\mbox{\rm d}}

\newtheorem{theore}{Theorem}
\renewcommand\thetheore{\arabic{section}.\arabic{theore}}
\newtheorem{theor}[theore]{\bf Theorem}
\newtheorem{rem}[theore]{Remark}
\newtheorem{propo}[theore]{\rm PROPOSITION}
\newtheorem{lem}[theore]{Lemma}
\newtheorem{definit}[theore]{\rm DEFINITION}
\newtheorem{coro}[theore]{\rm COROLLARY}
\newtheorem{exampl}[theore]{Example}
\newtheorem{examps}[theore]{Examples}
\newtheorem{case}{Case}

\renewcommand\thecase{({\it \roman{case}})}

\def\A{\mbox{${\cal A}$}}
\def\ga{\mbox{$\Delta(\cal A)$}}
\def\MA{\mbox{$M({\cal A})$}}
\def\gma{\mbox{$\Delta(M({\cal A}))$}}
\def\dg{\mbox{${\widehat{G}}$}}
\def\htg{\mbox{${\widetilde{H}(G)}$}}
\def\muhat{\mbox{$\widehat{\mu}$}}

\def\ccg{\mbox{${C_c(G)}$}}
\def\ckg{\mbox{${C_K(G)}$}}
\def\czg{\mbox{${C_o(G)}$}}
\def\cbg{\mbox{${C_b(G)}$}}
\def\hg{\mbox{${H(G)}$}}
\def\gdg{\mbox{${H(G)}$}}
\def\mg{\mbox{${M(G)}$}}
\def\czgw{\mbox{${C_o(G, 1/{\omega})}$}}
\def\cbgw{\mbox{${C_b(G, 1/{\omega})}$}}
\def\llocg{\mbox{${L_{\rm loc}(G)}$}}
\def\mlocg{\mbox{${M_{\rm loc}(G)}$}}
\def\l1gw{\mbox{${L^1(G, \omega)}$}}
\def\mgw{\mbox{${M(G, \omega)}$}}
\def\mzgw{\mbox{${M_0(G, \omega)}$}}
\def\mzzgw{\mbox{${M_{00}(G, \omega)}$}}
\def\hgw{\mbox{${H(G, \omega)}$}}
\def\glgw{\mbox{${\Delta(L^1(G, \omega))}$}}

\def\dl{\mbox{$\parallel$}}
\def\dlw{\mbox{$\parallel_{\omega}$}}
\def\dln{\mbox{$\parallel \cdot \parallel$}}
\def\dlnw{\mbox{$\parallel \cdot \parallel_{\omega}$}}

\title{A note on generalized characters}

\markboth{S J Bhatt and H V Dedania}{A note on generalized
characters}

\author{S J BHATT and H V DEDANIA}

\address{Department of Mathematics, Sardar Patel University,
Vallabh Vidyanagar~388~120, India\\
\noindent E-mail: subhashbhaib@yahoo.co.in; hvdedania@yahoo.com}

\volume{115}

\mon{November}

\parts{4}

\pubyear{2005}

\Date{MS received 23 November 2004; revised 9 July 2005}

\begin{abstract}
For a compactly generated LCA group $G$, it is shown that the set
$H(G)$ of all generalized characters on $G$ equipped with the
compact-open topology is a LCA group and $H(G) = \dg$ (the dual
group of $G$) if and only if $G$ is compact. Both results fail for
arbitrary LCA groups. Further, if $G$ is second countable, then
the Gel'fand space of the commutative convolution algebra $\ccg$
equipped with the inductive limit topology is topologically
homeomorphic to $\hg$.
\end{abstract}


\maketitle

\section{Introduction}

Throughout, let $G$ be a LCA group with Haar measure $\lambda$ and
let $\dg$ denote the dual group of $G$, i.e., the set of all
characters on $G$. Then it is well-known that $\dg$ is a LCA group
in compact-open topology. A {\it generalized character} on $G$ is
a continuous function $\alpha\hbox{\rm :}\  G \longrightarrow
{\mathbb C}^{\bullet}$, where \hbox{${\mathbb C}^{\bullet} = {\mathbb C}\!\setminus\!\{0\}$}
such that $\alpha(s+t) = \alpha(s)\alpha(t), s, t
\in G$. Let $H(G)$ denote the set of all generalized characters on
$G$ equipped with the compact-open topology. For $\alpha, \beta
\in H(G)$, define $(\alpha + \beta)(s) = \alpha(s)\beta(s), s \in
G$. Then $(H(G), +)$ is an abelian topological group~(23.34(b) of
\cite{HR}). It is straightforward to verify that $H({\mathbb Z})
\cong ({\mathbb C}^{\bullet}, \times)$ and $H({\mathbb T}) \cong
({\mathbb Z}, +)$, where ${\mathbb T}$ is the unit circle in
${\mathbb C}$.

Let $\ccg$ denote the set of all complex-valued continuous
functions on $G$ with compact support. Then $\ccg$ is a
commutative algebra with respect to the usual convolution product.
Let $\tau$ denote the inductive limit topology on $\ccg$. Then,
by~Lemma~2.1, p.~114 of \cite{M}, $(\ccg, \tau)$ is a commutative
topological algebra.

In this paper our main goal is to show that if $G$ is compactly
generated, then $\hg$ is a LCA group and that $H(G) = \dg$ if and
only if $G$ is compact. Both results fail for LCA groups. The
results appear to be a mathematical folklore; however we failed to
find a proof in the literature. In fact, the present note arises
out of our investigations of uniform norms in Beurling algebras
and weighted measure algebras~\cite{BhDe1,BhDe2}. As an
application we show that if, further, $G$ is second countable,
then the Gel'fand space $\Delta(\ccg)$ of $\ccg$ is homeomorphic
to $\hg$; in particular, $\Delta(\ccg)$ is a locally compact
space.

\section{Generalized characters}

\begin{lem}
\label{lem:tech} Let $m > 1$ be an integer and let $0 <
\varepsilon < 1/m$. Then there exists a natural number $N$ such
that{\rm ,} for each complex number $z$ satisfying $\varepsilon
\leq |z-1| \leq 1/m${\rm ,} there exists $1 \leq k \leq N$  such
that $|z^k-1| > 1/m$.
\end{lem}

\begin{proof}
For $r > 0$ and for $z \in {\mathbb C}$, let $\Gamma(z, r)$ denote
the circle with radius $r$ and center $z$. For $\delta > 0$, let
$L_{\delta} := \{r\hbox{e}^{i\delta}\hbox{\rm :}\ r > 0 \}$, the
open ray with angle $\delta$. Choose $0 < \delta < \pi/2$ such
that $L_{\delta}$ cuts the circle $\Gamma(1, \varepsilon)$ in two
points $z_0 = r_0\hbox{e}^{i\delta}$ and  $z_1 =
r_1\hbox{e}^{i\delta}$, where $r_0 < 1 < r_1$.

Now fix $z = r \hbox{e}^{i\theta}$ such that $\varepsilon \leq
|z-1| \leq 1/m$. Then $|\theta| < \pi/2$. Without loss of
generality, we may assume that $\theta \geq 0$. Then we have the
following three possibilities:

\begin{case}{\rm
$\theta \geq \delta$. Choose $n_1 \in {\mathbb N}$ such that
$n_1\delta \leq \pi/2$ and $L_{n_1\delta}$ does not intersect the
circle $\Gamma(1, 1/m)$. Then there exists $1 \leq k \leq n_1$
such that $|z^k-1| > 1/m$.}
\end{case}

\begin{case}{\rm
$r < r_{0}$. Choose $n_2 \in {\mathbb N}$ such that $r_0^{n_2}
< 1- 1/m$. Then $|z^{n_2}-1| \geq 1 - |z|^{n_2} = 1 - r^{n_2} > 1
- r_0^{n_2} \geq 1/m$.}
\end{case}

\begin{case}{\rm
$r_{1} < r$. Choose $n_3 \in {\mathbb N}$ such that $r_1^{n_3}
\geq 1+ 1/m$. Then $|z^{n_3} - 1| \geq |z|^{n_3} - 1 = r^{n_3} - 1
> r_1^{n_3} - 1 \geq 1/m$.}
\end{case}

Finally take $N = \max \{n_1, n_2, n_3 \}$. Then $N$ has the
required property.\hfill $\Box$
\end{proof}

\begin{theor}[\!]
\label{thm:gclca} Let $G$ be a compactly generated LCA group. Then
\begin{enumerate}
\renewcommand\labelenumi{{\rm (\roman{enumi})}}
\leftskip .15pc
\item $H(G)$ is a LCA group.
\item $H(G) = \dg$ if and only if $G$ is compact.
\end{enumerate}
\end{theor}

\begin{proof}$\left.\right.$
\begin{enumerate}
\renewcommand\labelenumi{(\roman{enumi})}
\leftskip .15pc
\item Fix an integer $m >1$. Define $V_m := \{z \in {\mathbb C}
\hbox{\rm :}\  |z-1| < 1/m \}$. Since $G$ is a compactly generated
LCA group, there exists a neighbourhood $U$ of $0$ in $G$ such
that its closure $\overline{U}$ is compact and it generates $G$
due to Theorem~5.13 of \cite{HR}. Take $T_m := N(\overline{U},
V_m) := \{\alpha \in \gdg\hbox{\rm :}\ \alpha(\overline{U})
\subseteq V_m \}$. Then $T_m$ is a neighbourhood of the identity
$1_G$ in $\gdg$. First we show that $T_m$ is equicontinuous at $0$
in $G$. Let $\varepsilon > 0$. If $\varepsilon \geq 1/m$, then $V
:= U$ is a neighbourhood of $0$ in $G$ such that
\begin{equation*}
s \in V\ \hbox{and}\ \alpha \in T_m \Longrightarrow
|\alpha(s)-\alpha(0)| = |\alpha(s)- 1| < 1/m \leq \varepsilon.
\end{equation*}
So we may assume that $\varepsilon < 1/m$. Then, by
Lemma~\ref{lem:tech}, one can find an integer $N$ such that, for
each $\varepsilon \leq |z-1| \leq 1/m$, there exists $1 \leq k
\leq N$  such that $|z^k-1| > 1/m$. Choose a neighbourhood $W$ of
$0$ in $G$ such that $\sum_{k=1}^{N}W_k \subseteq U$, where $W_k =
W$. Suppose, if possible, there exist $t \in W$ and $\alpha \in
T_m$ such that $|\alpha(t)-1| \geq \varepsilon$. Then, by the
definition of $N$, there exists $1 \leq k \leq N$ such that
$|\alpha(kt)-1| = |\alpha(t)^k-1| > 1/m$. On the other hand, $kt
\in U$ and so $|\alpha(kt)-1| \leq 1/m$. This is a contradiction.
Hence, we have
\begin{equation*}
s \in W\ \hbox{and}\ \alpha \in T_m \Longrightarrow |\alpha(s)-1|
< \varepsilon.
\end{equation*}
This proves that $T_m$ is equicontinuous at $0$ in $G$. Finally,
let $t \in G$ be arbitrary. Since $G$ is generated by $U$, there
exist $t_1, \dots, t_p \in U$ such that $t = t_1 + \cdots + t_p$.
Then, for each $\alpha \in T_m$,
\begin{equation*}
|\alpha(t)|=|\alpha(t_1)|\dots |\alpha(t_p)| \leq (1+1/m)^p.
\end{equation*}
By the above argument, one can choose a neighbourhood $W$ of $0$
in $G$ such that
\begin{equation*}
s \in W\ \hbox{and}\ \alpha \in T_m \Longrightarrow
|\alpha(s)-\alpha(0)| < \frac{\varepsilon}{(1+1/m)^p}.
\end{equation*}
Hence
\begin{align*}
|\alpha(s+t) - \alpha(t)| = |\alpha(s) - \alpha(0)||\alpha(t)|
\leq |\alpha(s) - 1|(1+1/m)^p < \varepsilon.
\end{align*}
\looseness 1 This proves that $T_m$ is equicontinuous. So its
closure $\hbox{Cl}_p(T_m)$ in the pointwise topology is
equicontinuous~(p.~17 of \cite{H}). Let $\hbox{Cl}_c(T_m)$ denote
the closure of $T_m$ in the compact-open topology. Then
$\hbox{Cl}_c(T_m) \subseteq \hbox{Cl}_p(T_m)$. Hence
$\hbox{Cl}_c(T_m)$ is equicontinuous.

\parindent 1pc Now take $t \in G$. Then $t = t_1 + \cdots + t_p$ for some $t_1,
\dots, t_p \in U$. Then $|\alpha(t)| =
|\alpha(t_1)|\cdots|\alpha(t_p)| \leq (1 + |\alpha(t_1)-1|) \cdots
(1 + |\alpha(t_p) - 1|) \leq (1+1/m)^p$ for each $\alpha \in T_m$.
Similarly, $|\alpha(t)| = |\alpha(t_1)|\cdots|\alpha(t_p)| \geq (1
- |\alpha(t_1)-1|) \cdots (1 - |\alpha(t_p) - 1|) \geq (1-1/m)^p$
for each $\alpha \in T_m$. Hence the closure of the set $T_m(t) :=
\{\alpha(t) \hbox{\rm :}\  \alpha \in T_m \}$ is compact in
${\mathbb C}^{\bullet}$. So, by Ascoli's theorem,
$\hbox{Cl}_c(T_m)$ is compact. This proves that $\gdg$ is a LCA
group.

\item \looseness -1 Let $G$ be compact and let $\alpha \in H(G)$. Since $\alpha$ is a
continuous group homomorphism, $\alpha(G)$ is a compact subgroup
of $({\mathbb C}^{\bullet}, \times)$. Hence $\alpha(G)$ is
contained in the unit circle. So $\alpha \in \dg$. For the
converse, assume that $H(G) = \dg$ and $G$ is compactly generated.
Then, by Theorem 9.8 of \cite{HR}, $G$ is topologically isomorphic
to ${\mathbb R}^m \times {\mathbb Z}^n \times K$ for some
non-negative integers $m, n$ and some compact group $K$. Then $\dg
= H(G) \cong H({\mathbb R}^m) \oplus H({\mathbb Z}^n) \oplus H(K)$
due to 23.34(c) of \cite{HR}. This implies that we must have $m
= n = 0$. So $G = K$ is compact. $\hfill \Box$\vspace{-1pc}
\end{enumerate}
\end{proof}

\begin{rem}
{\rm The following is an alternative proof of
Theorem~\ref{thm:gclca}(i). By the structure theory, a compactly
generated LCA group $G$ is a direct product of ${\mathbb R}^n$,
${\mathbb Z}^m$, and a compact group. By 23.34(c) of \cite{HR},
$H(G_1 \times G_2)$ is canonically homeomorphic to $H(G_1) \times
H(G_2)$. So it is enough to show that $H(G)$ is locally compact
for $G = {\mathbb Z}$ and $G = {\mathbb R}$. It is easy to see for
$G = {\mathbb Z}$. Observe that every continuous homomorphism
$\psi \hbox{\rm :}\ {\mathbb R} \longrightarrow {\mathbb C}$ is
differentiable and satisfies $\psi^{\prime}(t) = \psi(0)\psi(t), t
\in {\mathbb R}$, and so $\psi(t) = \exp(zt)$ for a unique complex
number $z$. Thus the map $\Lambda \hbox{\rm :}\ H({\mathbb R})
\longrightarrow ({\mathbb C}, +)$ is a bijective map. For $0 <
\varepsilon < 1$, let $W_{n, \varepsilon} = \{z \hbox{\rm :}\
|\hbox{e}^{zx} - 1| < \varepsilon, x \in [-n, n]\}$. Then it is
easy to see\break that
\begin{equation*}
W_{n, \varepsilon} \subseteq \{\alpha + i \beta \hbox{\rm :}\
|\alpha| < (1/n)\log(1+\varepsilon), |\beta| < (\cos^{-1}u_{n,
\varepsilon})/n \},
\end{equation*}
\looseness 1 where $u_{n, \varepsilon} =
(\hbox{e}^{-2|\alpha|n}+1-\varepsilon^{2})/(2\hbox{e}^{|\alpha|n})$.
Thus the mapping $\Lambda$ is open. Now for $0 < \delta < 1$,
\begin{equation*}
\{\alpha + i \beta \hbox{\rm :}\  |\alpha| < \log(1+\delta/2)/n,
|\beta| < \delta/2 \} \subseteq W_{n, \varepsilon}.
\end{equation*}
So $\Lambda$ is continuous. This completes the proof.\hfill
$\Box$}
\end{rem}

\begin{examps} \label{exam:e1}
{\rm The following two examples show that the above theorem is not
true for arbitrary LCA groups.
\begin{enumerate}
\renewcommand\labelenumi{{\rm (\roman{enumi})}}
\leftskip .15pc
\item Let $G = \{ \overline{n} = (n_1, \dots, n_k, 0, 0, \dots )
\hbox{\rm :}\  k \in {\mathbb N} \textrm{ and } n_i \in {\mathbb
Z} \}$ with the co-ordinatewise addition and the discrete
topology. Then $H(G) \cong {\mathbb C}^{\bullet {\mathbb N}}$ with
the pointwise topology. Then $H(G)$ is not a LCA group.

\item Let $G$ be an infinite abelian group having all elements of finite
order and the topology being the discrete topology. Let $\alpha
\in H(G)$ and let $s \in G$. Then there exists a natural number
$n$ such that $ns = 0$ and so $\alpha(s)^n = \alpha(ns) =
\alpha(0) = 1$, i.e., $|\alpha(s)| = 1$. Hence $\alpha \in \dg$.
Thus $H(G) = \dg$ and $G$ is not compact.
\end{enumerate}}
\end{examps}

\section{Gel'fand space of $\pmb{\ccg}$}

For $f \in \ccg$ and $t \in G$, let $(\tau_tf)(s) = f(s-t), t \in
G$. We know that, for $f \in \ccg$, the map $\Lambda_f \hbox{\rm
:}\  G \longrightarrow (\ccg, \|\cdot\|_1); s \longmapsto \tau_sf$
is continuous, where $\|\cdot\|_1$ is the $L^1$-norm. We prove the
following:

\setcounter{theore}{0}
\begin{lem}
\label{lem:1} Let $G$ be second countable{\rm ,} and let $f \in
\ccg$. Then the map $\Lambda_f \hbox{\rm :}\  G \longrightarrow
(\ccg, \tau); s \longmapsto \tau_sf$ is continuous.
\end{lem}

\begin{proof}
Since $G$ is a second countable, LCA group, $G$ is metrizable. Let
$d$ be an invariant metric on $G$ inducing the topology on $G$. So
it is enough to show that whenever $s_n \longrightarrow s$ in $G$,
we have $\Lambda_f(s_n) \longrightarrow \Lambda_f(s)$ in $\ccg$.
First, assume that $s = 0$. Let $U$ be a symmetric neighbourhood
of $0$ in $G$ such that $s_n \in U \; (n \in {\mathbb N})$ and
$\overline{U}$ is compact. Let $K = \overline{U} + {\rm supp}f$.
Then $K$ is compact, and the supports of $\tau_{s_n}f$ and $f$ are
contained in $K$.

Let $\varepsilon > 0$. Since $f|_{K}$ is continuous and since $K$
is a compact metric space, $f \hbox{\rm :}\  K \longrightarrow
{\mathbb C}$ is uniformly continuous. Let $\delta > 0$ such that
\begin{equation*}
s, t \in K \textrm{ and } d(s, t) < \delta \; \Longrightarrow \;
|f(s) - f(t)| < \varepsilon .
\end{equation*}
Choose  $n_0 \in {\mathbb N}$ such that $d(s_n, 0) < \delta \; (n
\geq n_0)$. Finally, let $t \in K$ and let $n \geq n_0$.

\setcounter{case}{0}
\begin{case}{\rm
$t-s_n \in K$: This implies $d(t-s_n, t) = d(-s_n, 0) = d(s_n, 0)
< \delta$; and so $|f(t-s_n) - f(t)| < \varepsilon$.}
\end{case}

\begin{case}{\rm
$t-s_n \notin K$: This implies $t \notin {\rm supp}f$; because if
$t \in {\rm supp}f$, then $t-s_n \in {\rm supp}f + \overline{U} =
K$ which is not the case. Hence $f(t-s_n) = f(t) = 0$; and so
$|f(t-s_n) - f(t)| < \varepsilon$.}
\end{case}

Hence $|\Lambda_f(s_n)(t) - \Lambda_f(0)(t)| = |f(t-s_n)-f(t)| <
\varepsilon , t \in K, n \geq n_0$. Thus $\|\Lambda_f(s_n) -
\Lambda_f(0)\|_{K} < \varepsilon \; (n \geq n_0)$. Thus
$\Lambda_f(s_n)\longrightarrow \Lambda_f(0)$.

Now let $s_n \longrightarrow s$ in $G$. Then $s_n - s
\longrightarrow 0$ in $G$. But $\|\Lambda_f(s_n) -
\Lambda_f(s)\|_{K} = \|\Lambda_f(s_n-s) - \Lambda_f(0)\|_{K}$.
Hence $\Lambda_f(s_n)\longrightarrow \Lambda_f(s)$. \hfill $\Box$
\end{proof}

Let $\Delta(\ccg)$ denote the Gel'fand space of $\ccg$. For
$\alpha \in H(G)$, define $\varphi_{\alpha}(f) =
\int_Gf(s)\alpha(s)\hbox{d}\lambda(s), f \in \ccg$. Then
$\varphi_{\alpha} \in \Delta(\ccg)$.

\begin{theor}[\!] \label{thm:1}
Let $G$ be second countable. Let $T \hbox{\rm :}\  H(G)
\longrightarrow \Delta(\ccg)$ be defined as $T(\alpha) =
\varphi_{\alpha}$. Then $T$ is a bijective continuous map.
\end{theor}

\begin{proof}
The mapping $T$ is clearly one-to-one. To show that $T$ is onto,
let $\varphi \in \Delta(\ccg)$. Then, for all $s \in G$ and for
all $f \in \ccg$,
\begin{equation*}
\varphi(f)^2 = \varphi(f^2) = \varphi(\tau_sf \ast \tau_{-s}f) =
\varphi(\tau_sf) \varphi(\tau_{-s}f).
\end{equation*}
This implies that if $\varphi(f) \neq 0$, then $\varphi(\tau_sf)
\neq 0$ for all $s \in G$. Let $f \in \ccg$ such that $\varphi(f)
\neq 0$. Define $\alpha \hbox{\rm :}\  G \longrightarrow {\mathbb
C}^{\bullet}$ as
\begin{equation*}
\alpha(s) = \frac{\varphi(\tau_sf)}{\varphi(f)}.
\end{equation*}
Note that $\alpha$ does not depend on $f$; because if $g \in \ccg$
is another function such that $\varphi(g) \neq 0$, then
\begin{equation*}
\varphi(\tau_sf)\varphi(g) = \varphi(\tau_sf \ast g) = \varphi(f
\ast \tau_sg) = \varphi(f)\varphi(\tau_sg), \quad s \in G.
\end{equation*}
Now, for $s, t \in G$,
\begin{equation*}
\alpha(s + t) = \frac{\varphi(\tau_{s+t}f)}{\varphi(f)} =
\frac{\varphi(\tau_{s}(\tau_tf))}{\varphi(f)} =
\frac{\varphi(\tau_{s}(\tau_tf))}{\varphi(\tau_tf)}\frac{\varphi(\tau_{t}f)}{\varphi(f)}
= \alpha(s)\alpha(t).
\end{equation*}
Since $G$ is second countable, the mapping $G \longrightarrow \ccg
; \; s \longmapsto \tau_sf$ is continuous due to
Lemma~\ref{lem:1}. Hence $\alpha$ is continuous. Thus $\alpha \in
H(G)$. Let $\mu \in M_{\rm loc}(G)$ be the Radon measure
corresponding to $\varphi$ (p.~838 of \cite{D}). Then, for $g \in
\ccg$,
\begin{align*}
\varphi_{\alpha}(g) & = \int_G g(s)\alpha(s)\hbox{d}\lambda(s)\\[.4pc]
& = \frac{1}{\varphi(f)}\int_G
g(s)\varphi(\tau_sf)\hbox{d}\lambda(s)\\[.4pc]
& = \frac{1}{\varphi(f)}\int_G
g(s)\int_Gf(t-s)\hbox{d}\mu(t)\hbox{d}\lambda(s)\\[.4pc]
& = \frac{1}{\varphi(f)}\int_G (f \ast g)(t)\hbox{d}\mu(t)\\[.4pc]
& = \frac{1}{\varphi(f)}\varphi(f \ast g) = \varphi(g).
\end{align*}
Thus $\varphi = \varphi_{\alpha}$. Hence $T$ is bijective. Now it
is easy to show that $T$ is continuous. $\hfill \Box$ \end{proof}

\begin{definit}$\left.\right.$\vspace{.5pc}

{\rm  \noindent For $\alpha \in \hg$, $\varepsilon > 0$, and
$\{f_1, \dots, f_n\} \subseteq \ccg$, define
\begin{equation*}
B(\alpha ; \varepsilon ; f_1, \dots, f_n) = \{\beta \in \hg
\hbox{\rm :}\  |\widehat{f_i}(\beta ) - \widehat{f_i}(\alpha )| <
\varepsilon \; (1 \leq i \leq n) \},
\end{equation*}
where $\widehat{f}(\beta) = \varphi_{\beta}(f) = \int_G
f(s)\beta(s)\hbox{d}\lambda(s)$. Then the collection
\begin{equation*}
{\cal B} = \{B(\alpha; \varepsilon; f_1, \dots, f_n) \hbox{\rm
:}\  \alpha \in \hg, \varepsilon > 0, \; n \in {\mathbb N}, \; \{f_1, \dots,
f_n\}\!\subseteq\!\ccg \}
\end{equation*}
forms a basis for some topology on $\hg$. Let $\tau_g$ denote the
topology on $\hg$ generated by this basis. Then $\tau_g \subseteq
\tau_{co}$ on $\hg$. Let $\htg$ denote the $H(G)$ equipped with
the topology $\tau_g$. We say that $\htg = \hg$ if $\tau_{co} =
\tau_g$.}
\end{definit}

\begin{rem} \label{lem:2}
{\rm Let $r > 1$. Define $\omega(s) = \hbox{e}^{r|s|}, s \in
{\mathbb R}$. Then $\omega$ is a weight on ${\mathbb R}$ such that
$\Delta(L^1({\mathbb R}, \omega)) \cong \Pi_{-r,r} := \{x+iy \in
{\mathbb C} \hbox{\rm :}\  -r \leq x \leq r \}$ due to~Theorem
4.7.33, p.~533 of \cite{D}.}
\end{rem}

\begin{theor}[\!] \label{thm:2}
If $G$ is {\rm (i)} discrete{\rm ,} {\rm (ii)} compact or {\rm
(iii)} $G = {\mathbb R}${\rm ,} then $\htg = \hg$.
\end{theor}

\begin{proof}
In the first two cases, it is enough to prove that the point
evaluation map $e\hbox{\rm :}\  G \times \htg \longrightarrow
{\mathbb C}$ is continuous due to~Corollary~13.1.1, p.~281 of
\cite{W}.

\begin{enumerate}
\renewcommand\labelenumi{(\roman{enumi})}
\leftskip .4pc
\item Fix $(g_0, \alpha_0)$ in $G \times \htg$. Let $V$ be a
neighbourhood of $e(g_0, \alpha_0) = \alpha_0(g_0)$ in ${\mathbb
C}$. Then there exists $\varepsilon > 0$ such that
$S(\alpha_0(g_0), \varepsilon) \subseteq V$. Choose $U = \{g_0\}$
and $f = \delta_{g_0}$. Define $B = B(\alpha_0; \varepsilon; f)$.
Then $U \times B$ is a neighbourhood of $(g_0, \alpha_0)$ in $G
\times \htg$. Then, for $(g, \alpha) \in U \times B$,
\begin{equation*}
|\alpha(g) - \alpha_0(g_0)| = |\alpha(g_0) - \alpha_0(g_0)| =
|\widehat{f}(\alpha) - \widehat{f}(\alpha_0)| < \varepsilon.
\end{equation*}
Hence $e(g, \alpha) = \alpha(g) \in V$. Thus the map $e$ is
continuous.

\item Since $G$ is compact, $\hg = \dg$. Suppose $\{t_{\gamma}\} \subset
G$ and $\{\alpha_{\gamma}\} \subset \htg$ are nets that converge
to $t$ and $\alpha$, respectively. Let $f \in \ccg$ such that
$\widehat{f}(\alpha) \neq 0$. Choose $\gamma_0$ such that
\begin{equation*}
|\widehat{f}(\alpha) - \widehat{f}(\alpha_{\gamma})| <
\frac{|\widehat{f}(\alpha)|}{2},\quad \gamma \geq \gamma_0.
\end{equation*}
Hence $|\widehat{f}(\alpha)| - |\widehat{f}(\alpha_{\gamma})| <
\frac{|\widehat{f}(\alpha)|}{2}$; and so
$\widehat{f}(\alpha_{\gamma}) \neq 0, \gamma \geq \gamma_0$. It is
elementary that, for $s \in G$ and for $\beta \in \htg$, $\beta(s)
\widehat{f}(\beta ) = [\tau_{s}(f)]^{\wedge}(\beta )$. Hence
\begin{equation*}
\alpha(s)  = \frac{[\tau_{s}(f)]^{\wedge}(\alpha)}{\widehat{f}(\alpha )}, \quad s \in G
\end{equation*}
and
\begin{equation*}
\alpha_{\gamma}(s)  = \frac{[\tau_{s}(f)]^{\wedge}(\alpha_{\gamma}
)}{\widehat{f}(\alpha_{\gamma} )}, \quad s \in G; \; \gamma \geq
\gamma_0 .
\end{equation*}
Since $\widehat{f}(\alpha_{\gamma}) \longrightarrow
\widehat{f}(\alpha)$, it is enough to prove that
$[\tau_{t_{\gamma}}(f)]^{\wedge}(\alpha_{\gamma}) \longrightarrow
[\tau_{t}(f)]^{\wedge}(\alpha)$. But
\begin{align*}
|[\tau_{t_{\gamma}}(f)]^{\wedge}(\alpha_{\gamma}) - [\tau_{t}(f)]^{\wedge}(\alpha)| & \leq
 |[\tau_{t_{\gamma}}(f)]^{\wedge}(\alpha_{\gamma}) -
[\tau_{t}(f)]^{\wedge}(\alpha_{\gamma})| \\[.2pc]
&\quad\,+ |[\tau_{t}(f)]^{\wedge}(\alpha_{\gamma}) -
[\tau_{t}(f)]^{\wedge}(\alpha )|\\[.2pc]
& \leq   \|[\tau_{t_{\gamma}}(f)]^{\wedge} -
[\tau_{t}(f)]^{\wedge}\|_1 \\[.2pc]
&\quad\, + |[\tau_{t}(f)]^{\wedge}(\alpha_{\gamma}) -
[\tau_{t}(f)]^{\wedge}(\alpha )|.
\end{align*}
The right-hand side tends to $0$ as $\gamma \longrightarrow
\infty$. Hence the map $e$ is continuous.

\item Note that $H({\mathbb R}) \cong {\mathbb C}$ and $\tau_{co}$
is exactly the usual topology ${\cal U}$ on ${\mathbb C}$. So we
need to prove that $\tau_g = {\cal U}$. Let $S(z, \varepsilon)$ be
an open sphere in ${\mathbb C}$ and let $w \in S(z, \varepsilon)$.
Let $r > 1$ such that $S(z, \varepsilon) \subset \Pi_{-r,r}$. By
Remark~\ref{lem:2}, there exists a weight $\omega$ on ${\mathbb
R}$ such that $\Delta(L^1({\mathbb R}, \omega)) \cong \Pi_{-r,r}$.
Since $C_c({\mathbb R})$ is dense in $L^1({\mathbb R}, \omega)$,
$\Delta((C_c({\mathbb R}), \|\cdot\|_{\omega})) \cong \Pi_{-r,r}$.
So choose $g_1, \dots, g_n$ in $L^1({\mathbb R}, \omega)$ and
$\delta > 0$ such that $B(w; \delta ; g_1, \dots, g_n) \subseteq
S(z, \varepsilon)$. Choose $f_1, \dots, f_n$ in $\ccg$ such that
$\|f_i - g_i\|_{\omega} < \frac{\delta}{3} \; (1 \leq i \leq n)$.
Now let $u \in B(w; \frac{\delta}{3}; f_1, \dots, f_n)$. Then,
for $1 \leq i \leq n$,
\begin{align*}
|\widehat{g_i}(u) - \widehat{g_i}(w)| &\leq |\widehat{g_i}(u) -
\widehat{f_i}(u)| + |\widehat{f_i}(u) - \widehat{f_i}(w)| +
|\widehat{f_i}(w) - \widehat{g_i}(w)|\\[.2pc]
&\leq \|f_i - g_i\|_{\omega} + |\widehat{f_i}(u) -
\widehat{f_i}(w)| + \|f_i - g_i\|_{\omega}\\[.2pc]
&< 2 \frac{\delta}{3} + \frac{\delta}{3} = \delta .
\end{align*}
Hence $u \in B(w; \delta; g_1, \dots, g_n)$. Thus
$B(w; \frac{\delta}{3} ; f_1, \dots, f_n) \subseteq S(z,
\varepsilon)$. Since $w$ is arbitrary, $S(z, \varepsilon) \in
\tau_g$. Hence the two topologies are identical. $\hfill
\Box$\vspace{-2pc}
\end{enumerate}
\end{proof}

\begin{theor}[\!] \label{thm:3}
If $\widetilde{H}(G_i) = H(G_i), i = 1, 2${\rm ,} then
$\widetilde{H}(G_1 \oplus G_2) = H(G_1 \oplus G_2)$.
\end{theor}

\begin{proof}
Let $G = G_1 \oplus G_2$. It is enough to prove that the point
evaluation map $e \hbox{\rm :}\  G \times \widetilde{H}(G)
\longrightarrow {\mathbb C}$ is continuous. Let $s = s_1 \oplus
s_2 \in G$ and $\alpha \in \widetilde{H}(G)$. Since $H(G) \cong
H(G_1) \oplus H(G_2)$, there exist $\alpha_1 \in H(G_1)$ and
$\alpha_2 \in H(G_2)$ such that $\alpha = \alpha_1 \oplus
\alpha_2$. Let $V$ be a neighbourhood of $e(s, \alpha) =
\alpha_1(s_1)\alpha_2(s_2)$. Choose $\varepsilon > 0$ such that
$S(\alpha_1(s_1), \varepsilon) \cdot S(\alpha_2(s_2), \varepsilon)
\subseteq V$. Since $\widetilde{H}(G_i) = H(G_i), i = 1, 2$, there
exist basic neighbourhoods $W_1 = U_1 \times B(\alpha_1; \delta_1;
f_1, \dots, f_m)$ of $(s_1, \alpha_1)$ in $G_1 \times H(G_1)$
and $W_2 = U_2 \times B(\alpha_2; \delta_2; h_1, \dots, h_n)$ of
$(s_2, \alpha_2)$ in $G_2 \times H(G_2)$ such that
\begin{equation*}
(t, \beta) \in W_1 \Longrightarrow \beta(t) \in S(\alpha_1(s_1), \varepsilon);
\end{equation*}
and
\begin{equation*}
(t, \beta) \in W_2 \Longrightarrow \beta(t) \in S(\alpha_2(s_2), \varepsilon).
\end{equation*}
Take $W = U \times B$, where $U = (U_1 \oplus U_2)$ and $B =
B(\alpha_1; \delta_1; f_1, \dots, f_m) \oplus B(\alpha_2;
\delta_2; h_1, \dots, h_n)$. Let $(s, \beta) \in W$. Then $s =
s_1 \oplus s_2$ for some $s_i \in U_i , i = 1, 2$ and $\beta =
\beta_1 \oplus \beta_2$ for some $\beta_1 \in B(\alpha_1;
\delta_1; f_1, \dots, f_m)$ and $\beta_2 \in B(\alpha_2;
\delta_2; h_1, \dots, h_n)$. So $\beta(s) =
\beta_1(s_1)\beta_2(s_2)$. Now, for all $1 \leq i \leq m$,
\begin{align*}
|\widehat{h_1}(\alpha_2) | |\widehat{f_i}(\beta_1 ) -
\widehat{f_i}(\alpha_1 )| &= |(f_i
\times h_1)^{\wedge}(\beta_1 \oplus \alpha_2) - (f_i \times h_1)^{\wedge}(\alpha )|
\\[.4pc]
&< \delta \leq \delta_1 |\widehat{h_1}(\alpha_2 ) |.
\end{align*}
Hence $\beta_1 \in B(\alpha_1; \delta_1; f_1, \dots, f_m)$; and
so $\beta_1(s_1) \in S(\alpha_1(g_1), \varepsilon)$. Similarly, we
can show that $\beta_2(s_2) \in S(\alpha_2(g_2), \varepsilon)$.
Hence $e(s, \beta) = \beta(s) = \beta_1(s_1)\beta_2(s_2) \in
S(\alpha_1(g_1), \varepsilon) \cdot S(\alpha_2(g_2), \varepsilon)
\subseteq V$. Thus the map $e$ is continuous. $\hfill \Box$\vspace{.5pc}
\end{proof}

\begin{coro}\label{cor:1}$\left.\right.$\vspace{.5pc}

\noindent If $G$ is compactly generated{\rm ,} then $\htg = \hg$.
\end{coro}

\begin{proof}
Since $G$ is compactly generated, $G \cong {\mathbb R}^m \times
{\mathbb Z}^n \times K$, where $m$ and $n$ are non-negative
integers and $K$ is a compact group due to Theorem~9.8 of \cite{HR}. Now
the result follows from Theorems~\ref{thm:2} and \ref{thm:3}.
$\hfill \Box$\vspace{.5pc}
\end{proof}

\begin{coro} \label{cor:2}$\left.\right.$\vspace{.5pc}

\noindent If $G$ is second countable and compactly generated{\rm
,} then $\hg \cong \Delta(\ccg)${\rm ,} and hence $\Delta(\ccg)$
is locally compact.\vspace{.5pc}
\end{coro}

\begin{proof}
The topology $\tau_g$ on $\hg$ is nothing but the Gel'fand topology
on $\ccg$. So the result follows from Theorem~\ref{thm:1} and
Corollary~\ref{cor:1}. $\hfill \Box$\vspace{.5pc}
\end{proof}

\begin{theor}[\!]
If $G$ is discrete{\rm ,} then $\hg \cong \Delta(\ccg)$.
\end{theor}

\begin{proof}
Define $T \hbox{\rm :}\  \hg
\longrightarrow \Delta(\ccg)$ as in Theorem~\ref{thm:1}. Since $G$
is discrete, $T$ is a bijective continuous map as in the proof of
Theorem~\ref{thm:1}. Let $\{\varphi_{\gamma}\}$ be a net in
$\Delta(\ccg)$ such that $\varphi_{\gamma} \longrightarrow
\varphi$ in $\Delta(\ccg)$. Let $\alpha_{\gamma}, \alpha \in \hg$
such that $T(\alpha_{\gamma}) = \varphi_{\gamma}$ and $T(\alpha) =
\varphi$. Then, for each $s \in G$,
\begin{equation*}
\alpha_{\gamma}(s) = \varphi_{\gamma}(\delta_s) \longrightarrow \varphi(\delta_s) =
\alpha(s).
\end{equation*}
Since $G$ is discrete, $\alpha_{\gamma}
\longrightarrow \alpha$ in $\hg$. Hence the result is proved.
$\hfill \Box$
\end{proof}

\begin{rem}$\left.\right.$
{\rm\begin{enumerate}\leftskip .15pc
\renewcommand\labelenumi{(\roman{enumi})}
\item Let $G$ be as in Example~\ref{exam:e1}(i). Then $\Delta(\ccg)
\cong H(G)$ is not locally compact.

\item If the condition ``second countable'' in Lemma~\ref{lem:1} can
be dropped, then the same can be dropped from
Corollary~\ref{cor:2}; in this case, $\Delta(\ccg)$ is locally
compact for all compactly generated LCA groups.
\end{enumerate}}\vspace{-1.5pc}
\end{rem}

\section*{Acknowledgements}

The authors gratefully acknowledge the financial support from the
University Grants Commission as a major research project, letter
No. F.8-9/2004(SR). The authors are thankful to the referee for
suggesting an alternative proof of Theorem~2.2(i) and for
constructive suggestions leading to the present version of the
paper.

\end{document}